\documentclass[a4paper,12pt,amsart,frenchb]{article}

\usepackage{amsmath,amsbsy,amsfonts,amssymb}
\usepackage[french]{babel}
\newcommand{\ass}[2]{\vskip0.3cm\noindent
{\bf {#1}}. { \sl {#2}}\vskip0.3cm\noindent
}
  
\begin{document}
 
  \title{  Une variante d'un r\'esultat de Aizenbud, Gourevitch, Rallis et Schiffmann  }
\author{J.-L. Waldspurger}
\date{28 octobre 2009}
\maketitle

 Soient $F$ un corps local non archim\'edien de caract\'eristique nulle, $W$ un espace vectoriel sur $F$ de dimension finie $\geq 1$, $<.,.>$ une forme bilin\'eaire sym\'etrique et non d\'eg\'en\'er\'ee sur $W$, $W=V\oplus U$ une d\'ecomposition orthogonale o\`u $U$ est une droite. On note $M$, resp. $G$, le groupe orthogonal de $W$, resp. $V$, et $M^0$, resp. $G^0$, le sous-groupe sp\'ecial orthogonal (plus exactement, on note ainsi les groupes de points sur $F$ de ces groupes alg\'ebriques). Le groupe $G$ s'identifie au sous-groupe des \'el\'ements de $M$ qui fixent $U$ point par point. On veut prouver
 
 \ass{Th\'eor\`eme 1}{Soient $\pi$, resp. $\rho$, une repr\'esentation admissible irr\'eductible de $M^0$, resp. $G^0$. Alors $dim_{{\mathbb C}}(Hom_{G^0}(\pi_{\vert G^0},\rho))\leq 1$.}
 
 Dans l'article [AGRS], les auteurs d\'emontrent plusieurs r\'esultats de ce genre et en particulier l'analogue du th\'eor\`eme ci-dessus o\`u les groupes sp\'eciaux orthogonaux sont remplac\'es par les groupes orthogonaux. Le cas des groupes sp\'eciaux orthogonaux a une certaine importance, en particulier si l'on s'int\'eresse \`a la conjecture locale de Gross-Prasad. La preuve du th\'eor\`eme 1 que l'on pr\'esente ci-dessous est une simple variante de celle de [AGRS]. On n'en d\'etaillera que les parties qui diff\`erent sensiblement de celle-l\`a. Je remercie vivement G. Henniart pour m'avoir signal\'e le probl\`eme et pour des remarques pertinentes sur une premi\`ere version de l'article.
 
 \bigskip
 
 On note $\mathfrak{g}$ l'alg\`ebre de Lie de $G^0$. Posons $\tilde{G}=G\times \{\pm 1\}$. Ce groupe agit sur $G^0$, $\mathfrak{g}$ et $ V$ par
 $$(g,\epsilon)x=gx^{\epsilon}g^{-1},\,\,(g,\epsilon)X=  \epsilon gXg^{-1},\,\,(g,\epsilon)v=\epsilon gv,$$
 pour $(g,\epsilon)\in \tilde{G}$, $x\in G^0$,  $X\in \mathfrak{g}$, $v\in V$. On a $\tilde{G}\subset \tilde{M}$ et, par cette inclusion, $\tilde{G}$ agit sur $M^0$.  Posons $e(V)=[\frac{dim(V)+1}{2}]$, que l'on consid\`ere comme un \'el\'ement de ${\mathbb Z}/2{\mathbb Z}$. Notons $\bar{G}$ le sous-groupe des \'el\'ements $(g,\epsilon)\in \tilde{G}$ tels que $det(g)=\epsilon^{e(V)}$. Notons $\chi$ le caract\`ere $(g,\epsilon)\mapsto \epsilon$ de $\bar{G}$ et, pour tout espace ${\cal S}$ sur lequel $\bar{G}$ agit, notons ${\cal S}^{\bar{G},\chi}$ le sous-espace des \'el\'ements qui se transforment sous l'action de $\bar{G}$ selon le caract\`ere $\chi$.  Pour tout espace topologique raisonnable $X$, on note  ${\cal S}(X)$ l'espace des fonctions sur $X$ \`a valeurs complexes, localement constantes et \`a support compact,  et ${\cal S}'(X)$ son dual. On a

\ass{Th\'eor\`eme 1'}{L'espace ${\cal S}'(M^0)^{\bar{G},\chi}$ est nul.}

Prouvons que ce th\'eor\`eme entra\^{\i}ne le th\'eor\`eme 1. On fixe $g\in G$ tel que $det(g)=(-1)^{e(V)}$. On note $\sigma$ l'antiinvolution $x\mapsto gx^{-1}g^{-1}$ de $M^0$. Le th\'eor\`eme 1' implique que toute distribution sur $M^0$ invariante par conjugaison par $G^0$ est invariante par $\sigma$. Le corollaire 1.1 de [AGRS] s'applique: pour $\pi$ et $\rho$ comme dans le th\'eor\`eme 1, on a
$$dim_{{\mathbb C}}(Hom_{G^0}(\pi_{\vert G^0},\rho^*)dim_{{\mathbb C}}(Hom_{G^0}((\pi^*)_{\vert G^0},\rho))\leq 1.$$
Il reste \`a prouver que 
$$(1) \qquad Hom_{G^0}(\pi_{\vert G^0},\rho^*)\simeq Hom_{G^0}(\pi_{\vert G^0},\rho)\simeq Hom_{G^0}((\pi^*)_{\vert G^0},\rho).$$
Fixons $\delta\in G$ tel que $det(\delta)=-1$. On d\'efinit $\pi^{\delta}(x)=\pi(\delta x\delta^{-1})$. Si $dim(W)$ est impaire, on a $\pi^*\simeq \pi\simeq \pi^{\delta}$. Si $dim(W)$ est paire, on a $\pi^*\simeq \pi$ ou $\pi^*\simeq \pi^{\delta}$. De m\^eme en rempla\c{c}ant $W$ et $\pi$ par $V$ et $\rho$. On a aussi
$$Hom_{G^0}(\pi^{\delta}_{\vert G^0},\rho^{\delta})=Hom_{G^0}(\pi_{\vert  G^0},\rho).$$
La relation (1) s'ensuit. $\square$

Prouvons le th\'eor\`eme 1'. On a:

\ass{Proposition}{Supposons que ${\cal S}'(M^0\times W)^{\bar{M},\chi}=\{0\}$. Alors ${\cal S}'(M^0)^{\bar{G},\chi}=\{0\}$.}

Preuve. Comme dans [AGRS] prop. 4.1, on fixe $e\in U$ non nul. Par descente de Frobenius, l'hypoth\`ese implique que ${\cal S}'(M^0)^{\bar{M}_{e},\chi}=\{0\}$, o\`u $\bar{M}_{e}$ est le fixateur de $e$ dans $\bar{M}$. Ce fixateur est l'ensemble des $(m,\epsilon)\in \bar{M}$ tels que $m=g\oplus \epsilon$ conform\'ement \`a la d\'ecomposition $W=V\oplus U$, avec $g\in G$. En introduisant l'\'el\'ement central $\epsilon_{M}$ dans $M$, on a aussi $m=\epsilon_{M}(g\oplus 1)$, avec un autre $g$. La condition $(m,\epsilon)\in\bar{M}$ signifie que $det(g)=\epsilon^{dim(W)+e(W)}$. Or $dim(W)+e(W)\equiv e(V)\,\,mod\,\,2{\mathbb Z}$. Donc $\bar{M}_{e}$ est l'ensemble des $(\epsilon_{M},1)(g,\epsilon)$, avec $(g,\epsilon)\in \bar{G}$. Puisque $(\epsilon_{M},1)$ agit trivialement sur $M^0$, on en d\'eduit ${\cal S}'(M^0)^{\bar{G},\chi}=\{0\}$. $\square$

D\'esormais, on oublie $M$ et $W$ et on va prouver ${\cal S}'(G^0\times V)^{\bar{G},\chi}=\{0\}$. On prouve simultan\'ement ${\cal S}'(\mathfrak{g}\times V)^{\bar{G},\chi}=\{0\}$. On raisonne par r\'ecurrence sur $d=dim(V)$. On le v\'erifie imm\'ediatement pour $d=1$. Pour $d=2$, l'action de $\bar{G}$ sur $G^0$ est triviale. On doit montrer que ${\cal S}'(V)^{\bar{G},\chi}=\{0\}$. L'action respecte la forme quadratique. Sur l'ouvert ${\cal V}$ des \'el\'ements $v$ tels que $<v,v>\not=0$, les orbites pour $\bar{G}$ sont les m\^emes que pour $G^0$. On en d\'eduit ais\'ement ${\cal S}'({\cal V})^{\bar{G},\chi}=\{0\}$. Reste l'ensemble $\Gamma$ des \'el\'ements $v$ tels que $<v,v>=0$. Si la forme quadratique est anisotrope, il est r\'eduit au point $0$ et le r\'esultat est clair. Sinon, dans une base convenable, c'est l'ensemble des $(x,y)$ tels que $xy=0$. Sur $\Gamma\setminus \{0\}$, on voit qu'un \'el\'ement $T\in {\cal S}'(V)^{\bar{G},\chi}$ est forc\'ement multiple de la distribution
$$f\mapsto \int_{F^{\times}}f(x,0)d^*x-\int_{F^{\times}}f(0,y)d^*y,$$
avec des mesures de Haar pour la multiplication. Or cette distribution ne se prolonge pas \`a $V$ tout entier en une distribution $G^0$-invariante (c'est l'exemple donn\'e dans l'introduction de [AGRS]). Donc $T$ est nulle sur $\Gamma\setminus \{0\}$. Et, comme pr\'ec\'edemment, $T$ ne peut pas avoir pour support le seul point $0$. Donc $T=0$. D\'esormais, on suppose $d\geq3$, donc le centre $Z$ de $G^0$ a au plus deux \'el\'ements.

\ass{Lemme}{Soit $T\in {\cal S}'(G^0\times V)^{\bar{G},\chi}$, resp. $T\in {\cal S}'(\mathfrak{g}\times V)^{\bar{G},\chi}$. Alors le support de $T$ est contenu dans $Z{\cal U}\times V$, resp ${\cal N}\times V$  .}

Cf. [AGRS], lemme 4.1 dont on utilise les notations. Soit $a\in G^0$ semi-simple non central. On peut d\'ecomposer $V$ en somme orthogonale $V=V_{+}\oplus V_{-}\oplus \oplus_{i\in I}V_{i}$. Pour tout $i\in I$, on a une suite d'extensions $F_{i}/F_{\pm i}/F$, o\`u $[F_{i}:F_{\pm i}]=2$. On admet formellement le cas o\`u $F_{i}/F_{\pm i}$ est "d\'eploy\'ee", c'est-\`a-dire $F_{i}=F_{\pm i}\oplus F_{\pm i}$. L'espace $V_{i}$ a une structure de $F_{i}$-espace vectoriel et est muni d'une forme hermitienne non d\'eg\'en\'er\'ee. L'\'el\'ement $a$ agit par multiplication par $1$ sur $V_{+}$, $-1$ sur $V_{-}$ et $a_{i}$ sur $V_{i}$, o\`u $a_{i}$ est un \'el\'ement de $F_{i}^{\times}$ tel que $a_{i}\not=\pm 1$ et $Norm_{F_{i}/F_{\pm i}}(a_{i})=1$. La composante neutre du commutant de $a$ dans $G^0$ est $G_{+}^0\times G_{-}^0\times \prod_{i\in I}G_{i}$, o\`u les $G_{i}$ sont des groupes unitaires (en un sens appropri\'e si $F_{i}=F_{\pm i}\oplus F_{\pm i}$). Les \'el\'ements de $G^0$ de partie semi-simple $a$ sont produits de $a$ et d'\'el\'ements de ce groupe. L'argument de [AGRS] nous ram\`ene \`a construire un \'el\'ement $(g,-1)\in \bar{G}$, fixant $a$, et v\'erifiant la propri\'et\'e suivante.  Consid\`erons $T_{+}\in {\cal S}'(G_{+}^0\times V_{+})^{G_{+}^0}$, $T_{-}\in {\cal S}'(G_{-}^0\times V_{-})^{G_{-}^0}$ et, pour tout $i\in I$, $T_{i}\in {\cal S}'(G_{i}\times V_{i})^{G_{i}}$. Alors $T=T_{+}\otimes T_{-}\otimes \otimes_{i\in I}T_{i}$ est fixe par l'action de $(g,-1)$. On prend $g_{+}\in G_{+}$ tel que $(g_{+},-1)\in \bar{G}_{+}$, $g_{-}\in G_{-}$ tel que $(g_{-},-1)\in \bar{G}_{-}$ et, pour tout $i\in I$,  un automorphisme antilin\'eaire de $V_{i}$ (relativement \`a l'extension $F_{i}/F_{\pm i}$) pr\'eservant la forme hermitienne. On prend $g=g_{+}\times g_{-}\times \prod_{i\in I}g_{i}$. D'apr\`es les r\'esultats de [AGRS] pour les groupes unitaires ou lin\'eaires, et d'apr\`es l'hypoth\`ese de r\'ecurrence, l'\'el\'ement $(g,-1)$ fixe $T$. Il fixe $a$. Il suffit de v\'erifier que $(g,-1)$ appartient \`a $\bar{G}$. Pour tout $i\in I$, $dim_{F}(V_{i})$ est paire et $det(g_{i})=(-1)^{dim_{F}(V_{i})/2}$. Parce que $a\in G^0$, $dim(V_{-})$ est paire et $det(g_{-})=(-1)^{e(V_{-})}=(-1)^{dim(V_{-})/2}$. On a $det(g_{+})=(-1)^{e(V_{+})}$. Donc $det(g)=(-1)^{e}$, o\`u $e=e(V_{+})+\frac{dim(V)-dim(V_{+})}{2}$. Or $e\equiv e(V)\,\,mod\,\,2{\mathbb Z}$. $\square$

Tout \'el\'ement de $Z$ \'etant invariant par $\bar{G}$, on est ramen\'e aux distributions sur ${\cal U}\times V$, que l'on descend par l'application de Cayley aux distributions sur ${\cal N}\times V$. Cela nous ram\`ene au probl\`eme sur l'alg\`ebre de Lie.

\ass{Lemme}{Soit $T\in {\cal S}'(\mathfrak{g}\times V)^{\bar{G},\chi}$.  Alors le support de $T$ est contenu dans $\mathfrak{g}\times \Gamma$.}

Cf. [AGRS] prop. 4.2. On fixe $v\in V$ avec $<v,v>\not=0$, on va montrer que ${\cal S}'(\mathfrak{g})^{\bar{G}_{v},\chi}=\{0\}$, o\`u $\bar{G}_{v}$ est le fixateur de $v$. Comme dans la preuve de: th\'eor\`eme 1' implique th\'eor\`eme 1, $\bar{G}_{v}$ est l'ensemble des $(\epsilon_{G},1)(g_{1},\epsilon)$, avec $(g_{1},\epsilon)\in \bar{G}_{1}$, o\`u on a \'ecrit $V=V_{1}\oplus Fv$. De nouveau, l'\'el\'ement $(\epsilon_{G},1)$ agit trivialement sur $\mathfrak{g}$, on peut le supprimer. On a $\mathfrak{g}=\mathfrak{g}_{1}\oplus V_{1}$ et l'action de $\bar{G}_{1}$ par restriction de celle de $\bar{G}$ est "la bonne". $\square$

On utilise maintenant l'argument de [AGRS] lemme 3.3 et remarque avant le lemme 3.4. Munissons $\mathfrak{g}$ de la forme bilin\'eaire sym\'etrique $(X,Y)\mapsto trace(XY)$. On d\'efinit deux transformations de Fourier partielles $T\mapsto {\cal F}_{\mathfrak{g}}T$ et $T\mapsto {\cal F}_{V}T$, la premi\`ere sur la variable dans $\mathfrak{g}$ (relativement \`a la forme ci-dessus), la seconde sur la variable dans $V$, relativement \`a la forme $<.,.>$. Elles conservent l'espace ${\cal S}'(\mathfrak{g}\times V)^{\bar{G},\chi}$. On a aussi deux repr\'esentations de Weil du groupe m\'etaplectique $\tilde{SL}_{2}$. Soit $T\in {\cal S}'(\mathfrak{g}\times V)^{\bar{G},\chi}$. Alors $T$ est \`a support dans ${\cal N}$ et ${\cal F}_{\mathfrak{g}}T$ l'est aussi. Donc $T$ est invariante par la premi\`ere repr\'esentation de Weil. De m\^eme, $T$ est invariante par la seconde repr\'esentation de Weil. En dimension impaire, une repr\'esentation de Weil ne se descend pas au groupe $SL_{2}$, l'\'el\'ement du noyau de la projection de $\tilde{SL}_{2}$ sur $SL_{2}$ agit par multiplication par $-1$ et tout invariant est nul. Donc $T$ est nulle si $dim(\mathfrak{g})$ est impaire ou si $dim(V)$ est impaire. On a pos\'e $d=dim(V)$. On a $dim(\mathfrak{g})=d(d-1)/2$. Donc $T$ est nulle si $d$ est impaire ou si $d\equiv 2\,\,mod\,\,4{\mathbb Z}$. On suppose maintenant $d\equiv 0\,\,mod \,\,4{\mathbb Z}$. On utilise la preuve de [AGRS] paragraphe 5. Elle montre que le support de nos distributions est contenu dans l'ensemble des $(X,v)$, $X\in {\cal N}$ et $v\in Q(X)$. Elle montre aussi qu'il suffit de fixer $X$ nilpotent et, en notant $\bar{G}_{X}$ son fixateur, de prouver le lemme suivant. On note $T\mapsto \hat{T}$ la transformation de Fourier dans ${\cal S}'(V)$, similaire \`a $T\mapsto {\cal F}_{V}(T)$.

\ass{Lemme}{On suppose $d\equiv 0\,\,mod\,\,4{\mathbb Z}$. Soit $T\in {\cal S}'(V)^{\bar{G}_{X},\chi}$. Supposons $T$ et $\hat{T}$ \`a support dans $Q(X)$. Alors $T=0$.}

Pour d\'emontrer cette assertion, on doit se d\'ebarrasser de l'hypoth\`ese sur $d$. On d\'efinit le sous-groupe $\underline{G}\subset \tilde{G}$ form\'e des $(g,\epsilon)\in \tilde{G}$ tels que $det(g)=\epsilon^{d}$. Remarquons que $\underline{G}=\bar{G}$ si $d\equiv 0\,\,mod\,\,4{\mathbb Z}$. Le lemme se g\'en\'eralise sous la forme

\ass{Lemme}{  Soit $T\in {\cal S}'(V)^{\underline{G}_{X},\chi}$. Supposons $T$ et $\hat{T}$ \`a support dans $Q(X)$. Alors $T=0$.}

Preuve. Le lemme 5.3 de [AGRS] reste valable, en remarquant que, pour des d\'ecompositions $V=V_{1}\oplus V_{2}$, $X=X_{1}\oplus X_{2}$, si $(g_{1},-1)\in \underline{G}_{1,X_{1}}$ et $(g_{2},-1)\in \underline{G}_{2,X_{2}}$, alors $(g_{1}g_{2},-1)$ appartient \`a $\underline{G}_{X}$. Cela nous ram\`ene au cas o\`u le couple $(V,X)$ est de l'un des types suivants.

$1^{er}$ cas. $d$ est impair, $V$ est muni d'une base $(e_{i})_{i=1,...,d}$, avec $<e_{i},e_{j}>=\nu(-1)^{i}\delta_{i,d+1-j}$ o\`u $\nu$ est une constante non nulle,  $Xe_{i}=e_{i-1}$ pour $i\geq2$, $Xe_{1}=0$. On d\'ecompose $V=V_{1}\oplus V_{0}\oplus V_{2}$, o\`u $V_{1}$ est engendr\'e par les $e_{i}$ pour $i\leq(d-1)/2$, $V_{0}$ est la droite port\'ee par $e_{(d+1)/2}$ et $V_{2}$ est engendr\'e par les $e_{i}$ pour $i\geq(d+3)/2$. Introduisons l'application
$$\begin{array}{ccc}{\cal S}(V)&\to&{\cal S}(V_{0})\\ f&\mapsto&f_{0}\\ \end{array}$$
d\'efinie par 
$$f_{0}(v_{0})=\int_{V_{1}}f(v_{1}+v_{0})dv_{1}.$$
Dans [AGRS], les auteurs prouvent qu'il existe $R\in {\cal S}'(V_{0})$ telle que $T(f)=R(f_{0})$ pour tout $f\in {\cal S}(V)$. Il y a un \'el\'ement  $a\in G^0$, appartenant au tore diagonal, tel que $aXa^{-1}=-X$. Posons $g=-a$. On a $(g,-1)\in \underline{G}_{X}$. Faisons agir trivialement cet \'el\'ement sur $V_{0}$. L'application ci-dessus est alors \'equivariante pour l'action de $(g,-1)$. Donc $T$ est invariante par cet \'el\'ement.

$2^{\grave{e}me}$ cas. $d\equiv 0\,\,mod\,\,4{\mathbb Z}$, $V$ est muni d'une base $(e_{i})_{i=1,...,d/2}\cup (f_{i})_{i=1,...,d/2}$, chaque sous-famille engendrant des lagrangiens, on a $<e_{i},f_{j}>=(-1)^{i}\delta_{i,d/2+1-i}$, $Xe_{i}=e_{i-1}$, $Xf_{i}=f_{i-1}$ pour $i\geq2$ et $Xe_{1}=Xf_{1}=0$. On d\'ecompose $V=V_{1}\oplus V_{2}$ , o\`u $V_{1}$ est engendr\'e par les $e_{i}$ et $f_{i}$ pour $i\leq d/4$ et $V_{2}$ est engendr\'e par les $e_{i}$ et $f_{i}$ pour $i\geq d/4+1$. Dans [AGRS], les auteurs prouvent que $T$ est multiple de la distribution
$$f\mapsto \int_{V_{1}}f(v_{1})dv_{1}.$$
 Il y a un \'el\'ement $g\in G^0$, appartenant au tore diagonal, tel que $gXg^{-1}=-X$. L'\'el\'ement $(g,-1)$ appartient \`a $\underline{G}_{X}$ et fixe  la distribution ci-dessus. Donc $T$ est invariante par cet \'el\'ement. $\square$

{\bf Remarque.}  Dans [AGRS], les auteurs consid\`erent aussi le deuxi\`eme cas avec $d\equiv 2\,\,mod\,\,4{\mathbb Z}$. C'est inutile. Par un changement de base (remplacer $e_{i}$ et $f_{i}$ par $e_{i}+f_{i}$ et $e_{i}-f_{i}$), ce cas se ram\`ene \`a une somme orthogonale de deux couples du premier cas.

Cela ach\`eve la preuve du th\'eor\`eme 1'.

{\bf Remarque.} On a suppos\'e la caract\'eristique de $F$ nulle. Selon un travail r\'ecent de Henniart, la d\'emonstration doit s'\'etendre au cas o\`u la caract\'eristique de $F$ est diff\'erente de $2$.

\bigskip

{\bf R\'ef\'erence}

[AGRS] A. Aizenbud, D. Gourevitch, S. Rallis, G. Schiffmann: {\it Multiplicity one theorems}, arXiv 07094215v1, mathRT

\bigskip

Institut de math\'ematiques de Jussieu-CNRS

175, rue du Chevaleret

75013 Paris

e-mail: waldspur@math.jussieu.fr

 \end{document}